\documentclass[12pt]{article}

\usepackage{amstext}
\usepackage{amsthm}
\usepackage{amsmath}
\usepackage{amssymb}
\usepackage{latexsym}
\usepackage{amsfonts}
\usepackage{graphicx}
\usepackage{texdraw}
\usepackage{graphpap}
\usepackage{color}
\usepackage[inline,nomargin]{fixme}

\usepackage[pagebackref,hypertexnames=true, colorlinks, citecolor=black, linkcolor=blue, urlcolor=red]{hyperref}

\usepackage[backrefs]{amsrefs}

\usepackage{latexsym}
\usepackage{amssymb}
\usepackage{euscript}

\let\cal=\mathcal      

\def\mcc{M\raise.5ex\hbox{c}C}
\def\mccarthy{M\raise.5ex\hbox{c}Carthy}

\def\eg{{\it e.g. }}
\def\ie{{\it i.e. }}

\def\h{{\cal H}}

\def\K{{\cal K}}
\def\M{{\cal M}}
\def\N{{\cal N}}

\def\vare{\varepsilon}

\let\i=\infty

\def\={\ = \ }    
\def\ot{\otimes}

\def\C{\mathbb C}

\def\D{\mathbb D}

\def\be{\setcounter{equation}{\value{theorem}} \begin{equation}}
\def\ee{\end{equation} \addtocounter{theorem}{1}}
\def\beq{\begin{eqnarray*}}
\def\eeq{\end{eqnarray*}}

\def\vs{\vskip 5pt}

\def\bp{{\sc Proof: }}
\def\ep{{}{\hfill $\Box$} \vskip 5pt \par}

\def\bl{\begin{lemma}}
\def\el{\end{lemma}}
\def\bt{\begin{theorem}}
\def\et{\end{theorem}}
\def\bprop{\begin{prop}}
\def\eprop{\end{prop}}
\def\bd{\begin{definition}}
\def\ed{\end{definition}}
\def\br{\begin{remark}}
\def\er{\end{remark}}
\def\bexer{\begin{exercise}}
\def\eexer{\end{exercise}}

\newtheorem{theorem}{Theorem}[section]
\newtheorem{prop}[theorem]{Proposition}
\newtheorem{lemma}[theorem]{Lemma}
\newtheorem{cor}[theorem]{Corollary}

\newtheorem{definition}[theorem]{Definition}

\def\gdel{G_\delta}

\def\norm#1{\| #1 \|}

\def\L{{\mathcal L}}

\def\gdel{G_\delta}

\def\vare{\varepsilon}

\def\M{{\mathbb M}}
\def\bh{{B(\h)}}
\def\bhd{\bh^d}
\def\bho{{B_1(\h)}}

\renewcommand\N{{\mathbb N}}
\def\mn{\M_n}
\def\mnd{\mn^d}
\def\Mnd{{\mathcal M}_n^d}

\def\md{{\mathbb M}^{[d]}}
\def\Md{{\mathcal M}^{[d]}}
\def\gds{G_\delta^\sharp}

\def\d{\delta}
\def\norm#1{\| #1 \|}

\numberwithin{equation}{section}
\title{Non-commutative Holomorphic  Functions on Operator Domains
\footnote{MSC  46H30, 32A70 Key Words: NC functions, free holomorphic functions, intertwining preserving }}
\author{Jim Agler
\thanks{Partially supported by National Science Foundation Grant
DMS 1361720}
\\ Dept. of Mathematics \\
U.C. San Diego\\ 
9500 Gilman Drive
\\
La Jolla, CA 92093
\and
John E. M\raise.5ex\hbox{c}Carthy
\thanks{Partially supported by National Science Foundation Grant  
DMS 1300280
}
\\ 
 Dept. of Mathematics \\
Washington University\\
 1 Brookings Drive \\
 St. Louis, MO 63130
}

\begin{document}

\bibliographystyle{plain}
\maketitle

Abstract: We characterize functions of $d$-tuples of bounded operators on a Hilbert space that 
are uniformly approximable by free polynomials on balanced open sets.

\section{Introduction}

What are non-commutative holomorphic functions in $d$ variables? 
The question has been studied since the pioneering work of  J.L. Taylor \cite{tay73}, but 
there is still no definitive answer.
The class should certainly contain the non-commutative, also called free, polynomials, 
\ie polynomials defined on $d$ non-commuting variables. It should be some sort of generalization of
the free polynomials, analogous to how holomorphic functions are generalizations of polynomials in commuting variables.
Just as in the commutative case, the class will depend on the choice of domain. In this note, we shall consider domains that are sets of $d$-tuples of operators on a Hilbert space.

One approach is to study non-commutative convergent power series on domains in $\bhd$ (where $\bhd$
means $d$-tuples
of bounded operators on some Hilbert space $\h$). This has been done systematically in 
G. Popescu's monograph \cite{po10b}, following on earlier work such as \cites{akv06,bgm06,po06,po08,po10}. 

Working with non-commutative power series is natural and appealing, but does present some difficulties.
One is that assuming {\em a priori} that the series converges uniformly is a strong assumption, and could be hard
to verify if the function is presented in some other form. On every infinite dimensional Banach space there is an entire
holomorphic function with finite radius of uniform convergence \cite[p. 461]{din99}. A second difficulty is dealing with domains that are not the domains of convergence of power series.

Another approach to non-commutative functions is the theory of nc-functions. Let $\mn$ denote the $n$-by-$n$ complex matrices,
which we shall think of as operators on a finite dimensional Hilbert space, 
and let $\md = \cup_{n=1}^\i \mnd$. 

If $x = (x^1, \dots, x^d) $ and $y = (y^1, \dots, y^d)$ are $d$-tuples of operators on the spaces
$\h$ and $\K$ respectively, we let $x \oplus y$ denote
the $d$-tuple $ (x^1 \oplus y^1, \dots, x^d \oplus y^d)$ on $\h \oplus \K$; and if $s \in B(\h,\K)$ and
$t \in B(\K,\h)$ we let
$sx$  and $xt$ denote respectively $(sx^1, \dots, sx^d)$ and $(x^1 t, \dots, x^d t)$.

\bd
\label{defa1}
A function $f$ defined on some set $D \subseteq \md$ is called
 graded if, for each $n$, $f$ maps $D \cap \mnd$ into $\mn$.

We say $f$ is an nc-function if it is graded and if, whenever $x, y \in D$ and
there exists a matrix $s$ such that $ s x = y s$, then 
$s f(x) = f(y) s$.
\ed

The theory of nc-functions has recently become a very active 
area of research, see \eg 
\cites{kvv09,hkms09,hkm11a, hkm11b, hm12, hkm12, pas14, ptd13}. D.  Kaliuzhnyi-Verbovetskyi  and V. Vinnikov have written a monograph
\cite{kvv14}  which develops the important ideas of the subject.

Nc-functions are {\em a priori} defined on matrices, not operators.
Certain formulas that represent them (such as \eqref{eqa4} below) can be naturally extended to operators.
This raises the question of how one can intrinsically characterize functions on $\bhd$ that are in some sense
extensions of nc-functions.

The purpose of this note is to show that on balanced domains in $\bhd$  there is an algebraic property -- intertwining preserving -- that together with an appropriate  continuity
is necessary and sufficient for a function to have a convergent power series, which in turn
is equivalent to the function being approximable by free polynomials on finite sets.
Moreover, it is a variation on the idea of an nc-function. On certain domains $\gds$ defined below, 
 the properties of intertwining preserving and continuity are equivalent in turn to the function being the unique extension of a bounded
nc-function.


\bd
\label{defa2}
Let $\h$ be an infinite dimensional  Hilbert space, let $D \subseteq \bhd$, and let $F: D \to \bh$.
We say that $F$ is intertwining preserving (IP) if:

(i) Whenever $x,y \in D$ and there exists
some bounded linear operator $T \in \bh$ such that $ T x = y T$, then $T F(x) = F(y) T$.

(ii) Whenever $( x_ n )$ is a bounded sequence  in $ D$,  and there exists
some invertible bounded linear operator $ s :  \h \to  \oplus \h $ such that
\[
s^{-1} \begin{bmatrix} x_1 & 0 & \cdots\\
0 & x_2& \cdots\\
\cdots&\cdots&\ddots\end{bmatrix} s \in D,
\] then
\[
F( s^{-1} \begin{bmatrix} x_1 & 0 & \cdots \\
0 & x_2 & \cdots \\
\cdots&\cdots&\ddots\end{bmatrix} s ) \=
s^{-1} \begin{bmatrix} F(x_1) & 0 & \cdots\\
0 & F(x_2) & \cdots \\
\cdots&\cdots&\ddots
  \end{bmatrix} s.
\]
\ed

Note that every free polynomial is IP, and therefore this condition must be inherited by any function that is a
limit of free polynomials on finite sets.
Nc-functions have the property that $f(x \oplus y) = f(x) \oplus f(y)$, and we would like 
to exploit the analogous condition (ii) of IP functions. To do this, we would like our domains to be closed
under direct sums. However, we can only do this by some identification of $\h \oplus \h$ with $\h$.

\bd
\label{defa4}
Let $\h$ be an infinite dimensional Hilbert space.
We say a set $D \subseteq \bhd$ is closed with respect to countable direct sums if, 
for every bounded sequence $x_1, x_2, \dots  \in D$, there is a unitary $u: \h \to \h \oplus \h \oplus \cdots $ 
such that the $d$-tuple $u^* (x_1 \oplus x_2 \oplus \cdots )u \in D$.
\ed

Two natural examples  are the sets
\be
\label{eqj55}
\{ x \in \bhd : \| x^1 \|, \dots, \| x^d \| < 1 \}
\ee
and
\be
\label{eqj6}
\{ x \in \bhd :  x^1 (x^1)^* + \dots +  x^d (x^d)^*  < I \}.
\ee

\bd
\label{defa3}
Let $F: \bhd \to \bh$. We say $F$ is sequentially strong operator continuous (SSOC) if,
whenever $x_n \to x$ in the strong operator topology on $\bhd$, then
$F(x_n)$ tends to $F(x)$ in the strong operator topology on $\bh$.
\ed

Since multiplication is sequentially strong operator continuous, it follows that every free
polynomial is SSOC, and  this  property is also inherited by  limits on  sets that are closed w.r.t. direct sums.

 Here is our first main result. Recall that a subset $B$ of a complex vector space is called balanced if whenever $x \in B$ and
$\alpha$ is in the closed unit disk $\overline{\D}$, then $\alpha x \in B$. 
\vs

{\bf Theorem \eqref{thmc1}.} {\em
Let $D$ be a balanced  open set in $\bhd$ that is closed with respect to countable direct sums, 
and let $F : D \to \bh$. The following are equivalent:

(i) The function $F$ is intertwining preserving and sequentially strong operator continuous.

(ii) There is a power series expansion
\[
\sum_{k=0}^\i P_k(x) 
\]
that converges absolutely 
 at each point $x \in D$ to $F(x)$, where each $P_k$ is a homogeneous free polynomial of
degree $k$.

(iii) The function $F$ is uniformly approximable on finite subsets of $D$ by free polynomials.
}

\vs

Let $\d$ be an $I \times J$ matrix
of free polynomials in $d$ variables, where $I$ and $J$ are any positive integers. Then
\be
\label{eqa1}
\gdel \ := \ 
\{ x \in \md : \| \d (x) \| < 1 \} ,
\ee
where if $x$ is a $d$-tuple of matrices acting on $\C^n$, then we calculate the norm of $\d(x)$
as  the operator norm from $(\C^n)^J$ to $(\C^n)^I$.
Notice that $G_{\delta_1} \cap G_{\delta_2} = G_{\d_1 \oplus \d_2}$, so those sets form a base
for a topology; we call this the free topology on $\md$.


For the rest of this paper, we shall fix $\h$ to be a separable infinite dimensional Hilbert space,
and let $\bho$ denote the unit ball in $\bh$. Let $\{e_1, e_2, \dots, \}$ be a fixed orthonormal basis of $\h$,
and let $P_n$ denote orthogonal projection onto $\vee \{ e_1, \dots, e_n \}$.

There is an obvious extension of \eqref{eqa1} to $\bhd$; we shall call this domain $\gds$.
\be
\label{eqa2}
\gds \ := \ 
\{ x \in \bhd : \| \d (x) \| < 1 \} .
\ee

 Both \eqref{eqj55}
and \eqref{eqj6} are of the form \eqref{eqa2} for an appropriate choice of $\delta$.
Note that every $\gds$ is closed with respect to countable direct sums.

By identifying $\mn$ with $P_n \bh P_n$, we can embed $\md$ in $\bhd$. 
If a function $F : \gds \to \bh$ satisfies $ F (x) = P_n F(x) P_n$ whenever $x = P_n x P_n$,
then $F$ naturally induces a graded function $F^\flat$ on $\gdel$.

Here is a slightly simplified version of Theorem~\eqref{thme1} (the assumption that $0$ goes to $0$ is unnecessary, but without it the statement is more complicated).

\vs

\bt
Assume that 
$\gds$ is connected and contains $0$. Then every bounded nc-function
 on
$\gdel$ that maps $0$ to $0$ has a unique extension to an SSOC IP function on $\gds$. The extension
has
a  series expansion in free polynomials that converges uniformly on $G_{t \d}^\sharp$ for each $t>1$.
\et

\section{Intertwining Preserving Functions}
\label{secb}

The normal definition of a an nc-function is a graded function $f$  defined on a set $D \subseteq \md$
such that $D$ is closed with respect to direct sums, and such that $f$ preserves direct sums and similarities,
\ie $f(x \oplus y) = f(x) \oplus f(y)$ and if $ x = s^{-1} y s$ then $f(x) = s^{-1} f(y) s$, whenever
$x,y \in D \cap \mnd$ and $s$ is an invertible matrix in $\mn$. The fact that on such sets $D$ this definition agrees with our earlier Definition~\eqref{defa1} is proved in \cite[Prop. 2.1]{kvv14}.

There is a subtle difference between the nc-property and IP, because of the r\^ole of $0$.
For an nc-function, $f(x \oplus 0) = f(x)\oplus 0$, but for an IP function, we have $f(x \oplus 0) = f(x) \oplus
f(0)$. If $f(0) = 0$, this presents no difficulty; but $0$ need not lie in the domain of $f$, and even if it does, 
it need not be mapped to $0$. 

Consider, for an illustration, the case $d=1$ and the function $f(x) = x+1$.
For each $n \in \N$, let $M_n$ be the $n$-by-$n$ matrix that is $1$ in the $(1,1)$ entry and $0$ elsewhere.
As an nc-function, we have
\[
f: 
\begin{bmatrix} 1 & 0  & \cdots &0\\
0 & 0 &  \cdots&0 \\
\vdots& \vdots  & \ddots & \vdots \\
0 & 0 &  \cdots & 0 \end{bmatrix}
\mapsto
\begin{bmatrix} 2 & 0  & \cdots& 0 \\
0 & 1  & \cdots & 0 \\
\vdots& \vdots & \ddots & \vdots \\
0 & 0  & \cdots  & 1 \end{bmatrix}
\]
But now, if we wish to extend $f$ to an IP function on $\bh$, what is the image
of the diagonal operator $T$ with first entry $1$ and the rest $0$?
We want to identify $T$ with $M_n \oplus 0$, and map it to $f(M_n) \oplus 0$ --- but then each
$n$ gives a different image.

In order to interface with the theory of nc-functions, we shall assume that all
our domains contain $0$. To avoid the technical difficulty we just described, we shall compose our functions 
with M\"obius maps to ensure that $0$ is mapped to $0$.


\bl
\label{lemb1}
If $F$ is an IP function on $D \subseteq \bhd$, and $P \in \bh$ is a projection, then, for all $c \in D$ satisfying 
$c = cP$ (or $c = Pc$) we have
\be
\label{eqb1}
a = PaP \ \Rightarrow \ F(a) \= P F(a) P  + P^\perp F(c) P^\perp.
\ee
\el
\bp
As $Pa = a P$, we get $P F(a) = F(a) P$.
As $P^\perp a = 0 = c P^\perp$, we get $P^\perp F(a) = F(c) P^\perp$.
Combining these, we get \eqref{eqb1}.
\ep

We let $\phi_\alpha$ denote the M\"obius map on $\D$ given by
\[
\phi_\alpha(\zeta) \= \frac{\zeta - \alpha}{1 - \bar \alpha \zeta} .
\]

\bl
\label{lemb2} Let $D \subseteq \bhd$ contain $0$, and assume $F$ is an IP function from $D$ to $\bho$.
Then:

(i) $F(0) = \alpha I_\h$.

(ii) The map $ H(x) := \phi_\alpha \circ F (x)$ is an IP function on $D$ that maps $0$ to $0$.

(iii) For any $a \in D$ and any projection $P$ we have
\[
a = PaP  \ \Rightarrow \  H(a) = P H(a) P  .
\]

(iv) $F = \phi_{-\alpha} \circ H$.
\el
\bp
(i) By Lemma~\eqref{lemb1} applied to $a = c = 0$, we get that $F(0)$ commutes with every projection $P$
in $\bh$. Therefore it must be a scalar.

(ii) For all $z$ in $\bho$, we have 
\be
\label{eqb2}
\phi_\alpha(z) \= -  \alpha I_\h + (1 - | \alpha |^2) \sum_{n=1}^\i \bar \alpha^{n-1} z^n ,
\ee
where the series converges uniformly and absolutely on every ball of radius less than one.
By (i), we have $H(0) = \phi_\alpha (\alpha I_\h) = 0$. If $T x = yT$, then $T F(x) = F(y)T$,
and so $ T [ F(x)]^n = [F(y)]^n T$ for every $n$.
Letting $z$ be $F(x)$ and $F(y)$ in \eqref{eqb2} and using the fact that the series converges uniformly,
we conclude that $T \phi_\alpha (F(x) ) =  \phi_\alpha (F(y) ) T$, and hence $H$ is IP.

Now (iii) follows from Lemma~\eqref{lemb1} with $ c = 0$. We get (iv) because $\phi_{-\alpha} \circ
\phi_{\alpha} (z) = z$ for every $z \in \bho$.
\ep
By choosing a basis $\{e_1, e_2, \dots, \}$ for $\h$, we can identify
$\mnd$ with $P_n \bhd P_n$. Let us define
\[
\Mnd \ := \ P_n \bhd P_n,\quad \Md = \cup_{n=1}^\i \Mnd .
\]
Applying Lemma~\eqref{lemb2}, we get the following.
\bprop
\label{propb1}
Let $D \subseteq \bhd$ contain $0$, and assume $F$ is an IP function from $D$ to $\bho$. Let $H = \phi_\alpha \circ F$, where $\alpha $ is the scalar such that $F(0) = \alpha I_\h$.
Then $ H |_{D \cap \Md}$ is an nc-function that is bounded by $1$ in norm,
and maps $0$ in $\Mnd$ to the matrix $0$ in
${\mathcal M}_n^1 = P_n \bh P_n$.
\eprop

If we let $H^\flat$ denote  $ H |_{D \cap \Md}$, we can ask 

\begin{quest}
 To what extent does $H^\flat$ determine $H$?
\end{quest}

\begin{quest}
 Does every bounded nc-function from ${D \cap \Md}$ to ${\mathcal M}^1$ extend to a bounded IP function
on $D$?
\end{quest}

\vs
If \[
\delta(x^1,x^2) \= I - (x^1 x^2 - x^2 x^1) ,
\]
then $\gds$ is non-empty, but $\gdel$ is empty, and the questions do not make much sense.
 But we do give  answers
to both questions in Theorem~\eqref{thme1}, in the special case that $D$ is of the form $ \gds$ and in addition  is assumed to be balanced.

\section{IP SSOC functions are analytic}

Let us give a quick summary of what it means for a function to be holomorphic on a Banach space;
we refer the reader to the book \cite{din99} by S. Dineen for a comprehensive treatment.
Let $D$ be an open subset of a Banach space $X$, and $f: D \to Y$ a map into a Banach space $Y$.
We say $f$ has a G\^ateaux derivative at $x$ if
\[
\lim_{\lambda  \to 0} \frac{f(x + \lambda  h) - f(x)}{\lambda } \ := \ Df(x)[h] 
\]
exists for all $h \in X$.
If $f$ has a G\^ateaux derivative at every point of $D$ it is  G\^ateaux holomorphic 
\cite[Lemma 3.3]{din99}, \ie holomorphic on each one dimensional slice. If in addition
$f$ is locally bounded on $D$, then it is actually  Fr\'echet holomorphic  \cite[Prop. 3.7]{din99}, 
which means that for each $x$ there is a neighborhood $G$ of $0$ such that  the Taylor series
\be
\label{eqj1}
 f(x+h) \= f(x) + \sum_{k=1}^\i D^k f(x) [h,\dots, h]  \quad \forall h \in G,
\ee
 converges uniformly for all $h$ in $G$. The $k^{\rm th}$ derivative is a continuous
linear map from $X^k \to Y$, which is evaluated on the $k$-tuple $(h,h, \dots, h)$.

The following lemma is 
the IP version of  \cite[Prop. 2.5]{hkm11b} and \cite[Prop 2.2]{kvv14}.
\bl
\label{lemk1}
Let $D$ be an open set
in $\bhd$ that is closed with respect to countable direct sums, 
and let $F : D \to \bh$ be
intertwining
preserving. Then $F$ is bounded on bounded subsets of $D$, continuous and G\^ateaux differentiable.
\el
\bp 
(Locally bounded)
Suppose there were $x_n \in D$ 
such that $\{ \| x_n \| \} $ is bounded, but $\{ \| F(x_n) \| \}$ is unbounded.
Since $D$ is closed with respect to countable direct sums, there exists some unitary $u : \h \to \h^\i$ such that
$u^* (\oplus x_n ) u \in D$. Since $F$ is IP, by Definition \eqref{defa2}, we have
$ [\oplus F(x_n)]$ is bounded, which is a contradiction.

(Continuity)
Fix $a \in D $ and let $\vare > 0$. 
By hypothesis, there exists a unitary $u : \h \to \h^2$ such that
\be
\label{eqs2}
\alpha \ := \ u^*  \begin{bmatrix} a & 0 \\
0 & a \end{bmatrix} u \ \in \ D.
\ee
Choose $\delta_1 > 0$ such that $B( a , \delta_1) \subseteq D$,
$B( \alpha , \delta_1) \subseteq D$, and such that on $B( \alpha, \delta_1)$ the function $F$ is bounded by $M$.
Choose $\delta_2 > 0$ such that $\delta_2 < \min ( \frac{1}{2} \delta_1,  \frac{\vare}{2M} \delta_1)$.
Note that for any $a,b \in \bhd$ and any $\lambda \in \C$, we have
\be
\label{eqs1}
u^* 
\begin{bmatrix} I & - \lambda \\
0 & I \end{bmatrix}
\begin{bmatrix} b & 0 \\
0 & a \end{bmatrix}
\begin{bmatrix} I &  \lambda \\
0 & I \end{bmatrix}
u \=
u^*
\begin{bmatrix} b & \lambda(b-a) \\
0 & a \end{bmatrix}
u
\ee
So by part (ii) of the definition of IP (Def. \eqref{defa2}) we get that if $\| b - a \| < \delta_2$,
and letting $\lambda = \frac{M}{\vare}$,
then
\[
F( u^* 
\begin{bmatrix} I & - \frac{M}{\vare} \\
0 & I \end{bmatrix}
\begin{bmatrix} b & 0 \\
0 & a \end{bmatrix}
\begin{bmatrix} I &  \frac{M}{\vare} \\
0 & I \end{bmatrix}
u ) \=
u^*
\begin{bmatrix}F( b) &\frac{M}{\vare}[F (b)-F(a)] \\
0 & F(a) \end{bmatrix}
u
\]
is bounded by $M$.
In particular, since the norm of the $(1,2)$-entry of the last matrix is bounded by the norm of the whole matrix,
 we see that $\norm{(M/\vare)(F(b)-F(a))} < M$, so
$\norm{F(b)-F(a)} < \vare$.

\vs
(Differentiability) Let $a \in D$ and $h \in \bhd$. Let $u$ be as in \eqref{eqs2}.
Choose $\vare > 0$ such that, for all
complex numbers $t$ with $|t| < \vare$, 
\[
 u^*  \begin{bmatrix} a + th & \vare h \\
0 & a \end{bmatrix} u \ \in \ D,
\]
and $a + th \in D$.
Let $b = a +th$ and $\lambda = \frac{\vare}{t}$ in \eqref{eqs1}, and as before we conclude that
\be
\label{eqs3}
F(u^*  \begin{bmatrix} a + th & \vare h \\
0 & a \end{bmatrix} u )
\=
u^*  \begin{bmatrix} F(a + th) & \vare\, \frac{F(a+th) - F(a)}{t} \\
0 & F(a) \end{bmatrix} u .
\ee
As $F$ is continuous, when we take the limit as $t \to 0$ in \eqref{eqs3}, we get
\[
F(u^*  \begin{bmatrix} a  & \vare h \\
0 & a \end{bmatrix} u )
\=
u^*  \begin{bmatrix} F(a ) & \vare DF (a)[h] \\
0 & F(a) \end{bmatrix} u .
\] 
Therefore $DF(a) [h]$ exists, so $F$ is G\^ateaux differentiable, as required.
\ep

When we replace $X$ by a Banach algebra (in our present case, this is $\bhd$ with coordinate-wise multiplication), 
we would like something more than Fr\'echet holomorphic: 
we would like the $k^{\rm th}$ term in \eqref{eqj1} to be an actual free polynomial, homogeneous
of degree $k$, in the entries of $h$. 
%
%

The following result was proved by 
Kaliuzhnyi-Verbovetskyi  and Vinnikov \cite[Thm. 6.1]{kvv14}
and by I. Klep and S. Spenko \cite[Prop. 3.1]{ks14}.

\bt
\label{thmj11}
Let 
\begin{align}
\nonumber
g: \md &\to \M^1 \\
x&\mapsto g(x)
\nonumber
\end{align}
 be an nc-function such that each matrix entry of $g(x)$  is a polynomial of degree less than or equal to $N$ 
in the entries of the matrices $x^r, 1 \leq r \leq d$.
Then $g$ is a free polynomial of degree less than or equal to $N$.
\et
We extend this result to multilinear SSOC IP maps. Each $h_j$ will be a $d$-tuple of operators,
$(h_j^1, \dots, h_j^d)$.
\bprop
\label{propj12}
Let
\begin{align}
\nonumber
L: \bh^{dN} &\to \bh \\
(h_1, \dots, h_N) &\mapsto L(h_1,\dots, h_N) 
\nonumber
\end{align}
be a continuous $N$-linear map from $(\bh^{d})^N$ to $\bh$ that is IP and SSOC.
Then $L$ is a homogeneous polynomial of degree $N$ in the variables
$h_1^1, \dots, h_N^d$.
\eprop
\bp
By Proposition~\eqref{propb1}, if we restrict $L$ to ${\mathcal M}^{dN}$,
we get an nc-function. By Theorem~\eqref{thmj11}, there is a free polynomial $p$
of degree $N$ that agrees with $L$ on ${\mathcal M}^{dN}$.
By homogeneity, $p$ must be homogeneous of degree $N$.
Define
\[
\Delta(h) \= L(h) - p(h) .
\]
Then $\Delta$ 
 vanishes on $({\mathcal M}^{d})^{N}$, and is SSOC.
Since $({\mathcal M}^{d})^{N}$ is strong operator topology dense in $(\bh^{d})^N$,
it follows that 
 $\Delta$ is identically $0$.
\ep

%

One of the  achievements of Kaliuzhnyi-Verbovetskyi  and  Vinnikov in \cite{kvv14} is the Taylor-Taylor formula, \cite[Thm. 4.1]{kvv14}.
This comes with a remainder term, which can be estimated. They show \cite[Thm 7.4]{kvv14} that 
with the assumption of local boundedness, this renders an nc-function analytic.
The following theorem is an IP version of the latter result.

\bt
\label{thmj1}
Let $D$ be an open 
 neighborhood of $0$ in $\bhd$, and let $F : D \to \bh$ be a 
  function that is intertwining
preserving and sequentially strong operator continuous. 
Then there is an open set $U \subseteq D$ containing $0$
and homogeneous free polynomials $P_k$ of degree $k$ such that
\[
 F(x) \= F(0) + \sum_{k=1}^\i P_k(x)  \quad \forall x \in U,
\]
where the convergence is uniform for $x \in U$.
\et

\bp
Any open ball centered at $0$ is closed with respect to countable direct sums, 
so we can assume without loss of generality that $D$ is closed with respect to countable direct sums and bounded.
By Lemma~\eqref{lemk1}, $F$ is bounded and G\^ateaux differentiable on $D$, and  so by
\cite[Prop. 3.7]{din99}, 
 $F$ is automatically 
Fr\'echet holomorphic.
Therefore there is some open ball $U$ centered at $0$ such that
\[
F(h) \= F(0) + \sum_{k=1}^\i D^k F(0)[h,\dots, h]  \quad \forall h \in U .
\]
 We must show that each  $D^k F(0) [h,\dots, h]$ is actually a free polynomial in $h$.

Claim 1: For each $k \in \N$, the function 
\be
\label{eqj17}
G^k: (h^0, \dots, h^k) \ \mapsto \
 D^k F(h^0) [h^1,\dots, h^k]
\ee
is an IP
function on $ U \times (\bhd)^k  \subseteq (\bh^d)^{k+1}$.

Indeed, when $k=1$, we have
\be
\label{eqj5}
 D F(h^0) [h^1]  \= \lim_{t \to 0} \frac{1}{t} [ F(h^0 + th^1) - F(h^0)].
\ee
As $F$ is IP, so is the right-hand side of  \eqref{eqj5}.
For $k > 1$, 
\begin{align*}
 D^k F(h^0) [h^1, \dots, h^k]  \= \lim_{t \to 0} \frac{1}{t} &\big[D^{k-1} F(h^0 + th^k)[h^1, \dots, h^{k-1}]
\\
& - D^{k-1}F(h^0)[h^1, \dots, h^{k-1}] \big].
\end{align*}
By induction, these are all IP.

Claim 2:  For each $k \in \N$, the function 
$G^k$ from \eqref{eqj17} is SSOC 
on $ U \times (\bhd)^k$.

Again we do this by induction on $k$. Let $G^0 := F$, which is SSOC on $U \subseteq D$ by hypothesis.
 Since $G^{k-1}$ is IP on the  set
$U^k$,
 it is locally bounded, and by Lemma~\eqref{lemk1} it is G\^ateaux differentiable. 
Suppose 
\[
{\rm SOT} \lim_{n \to \i} h^j_n \ = \ h^j, \ 0 \leq j \leq k ,
\]
where each $h^j_n$ and $ h^j $ is in $U$.
Let $ h$ denote the $(k+1)$-tuple $ (h^0, \dots, h^{k})$ in $U^{k+1}$, and let
$\tilde h$ denote  the $k$-tuple $ (h^0, \dots, h^{k-1})$; similarly,
 let  $h_n$ denote $ (h_n^0, \dots, h_n^{k})$
 and $\tilde h_n$ denote
 $ (h_n^0, \dots, h_n^{k-1})$.
There exists some unitary $u$ so that
$ y = u^* (\tilde  h \oplus \tilde h_1 \oplus \tilde h_2 \oplus \cdots ) u$ is in $U^{k}$.
Since $G^{k-1}$ is differentiable at $y$, and is IP, we have that the diagonal operator with entries
\begin{align}
\nonumber
 &\frac{1}{t}
 [ G^{k-1}(h^0 + t h^k, h^1, \dots, h^{k-1}) - G^{k-1}(h^0 , h^1, \dots, h^{k-1})], \\
 & \frac{1}{t}
 [  G^{k-1}(h^0_1 + t h^k_1, h^1_1, \dots, h^{k-1}_1) - G^{k-1}(h^0_1 , h^1_1, \dots, h^{k-1}_1) ],
 \nonumber
 \\
 & \dots 
 \label{eqj19}
\end{align} 
has a limit as $t \to 0$.

Let $\vare > 0$, and let $v \in \h$ have $\| v \| \leq 1$.  Choose $t$ sufficiently close to $0$ that each of the difference quotients
in \eqref{eqj19} is within $\vare/3$ of its limit (which is $G^k$ evaluated at the appropriate $h$ or $h_n$).
Let $n$ be large enough so that 
\begin{align*}
&\| \big [ G^{k-1}(h^0 + t h^k, h^1, \dots, h^{k-1}) -  G^{k-1}(h^0_n + t h^k_n, h^1_n, \dots, h^{k-1}_n) \big] v \|  \\
&+ \ \| \big[ G^{k-1}(h^0 , h^1, \dots, h^{k-1}) - G^{k-1}(h^0_n , h^1_n, \dots, h^{k-1}_n) \big] v \|\\
 &\ \leq \ \frac{\vare t}{3} .
\end{align*}
Then
\[
\| \big[ G^k (h^0, \dots, h^k) - G^k (h^0_n, \dots, h^k_n ) \big] v \| \ \leq \ \vare .
\]
So each $G^k$ is SSOC on $U^{k+1}$. 
As $G^k$ is linear in the last $k$ variables, it is SSOC on $ U \times (\bhd)^k$
 as claimed.

Therefore for each $k$, the map 
\[
(h^1, \dots, h^k) \ \mapsto \
 D^k F(0) [h^1,\dots, h^k]
\]
is a linear IP
function that is SSOC in a neighborhood of $0$, so by Proposition~\eqref{propj12} is a free
polynomial.
\ep

\section{Power series}

\bt
\label{thmc1}
Let $D$ be a balanced  open set in $\bhd$ that is closed with respect to countable direct sums, 
and let $F : D \to \bh$. The following are equivalent:

(i) The function $F$ is intertwining preserving and sequentially strong operator continuous.

(ii) There is a power series expansion
\be
\label{eqc1}
\sum_{k=0}^\i P_k(x) 
\ee
that converges absolutely 
 at each point $x \in D$ to $F(x)$, where each $P_k$ is a homogeneous free polynomial of
degree $k$.

(iii) The function $F$ is  uniformly approximable on finite subsets of $D$ by free polynomials.

\et
\bp
(i) $\Rightarrow$ (ii). 
 As $F$ is bounded on bounded subsets of $D$ by Lemma~\eqref{lemk1}, it is Fr\'echet
  holomorphic. By Theorem~\eqref{thmj1}, the
power series at $0$ is actually of the form \eqref{eqc1}. We must show the series converges on all of $D$.

Fix $x \in D$. Since $D$ is open and balanced, there exists $r > 1$ such that $\lambda x \in D$
for every $\lambda \in \D(0, r)$.
As each $P_k$ is homogeneous, we have that for $\lambda$ in a neighborhood of $0$,
\be
\label{eqc2}
F(\lambda x) \= \sum_{k=0}^\i P_k( \lambda x) \= \sum_{k=0}^\i \lambda^k P_k(x) .
\ee
Therefore the function $\psi: \lambda \mapsto F(\lambda x)$ is analytic on $\D(0, r)$, with values in $\bh$,
and its power series expansion at $0$ is given by \eqref{eqc2}.
Let $M = \sup \{ \| F(\lambda x) \| : | \lambda | < r \}$.

By the Cauchy integral formula, since $\| F \|$ is bounded by $M$, we get that
\be
\label{eqc3}
\| \frac{d^k \psi}{d \lambda^k} (0) \| \ \leq \ M \frac{k!}{r^k} .
\ee
Comparing \eqref{eqc2} and \eqref{eqc3}, we conclude that
\[
\| P_k (x) \| \ \leq \ \frac{M}{r^k},
\]
and so the power series in \eqref{eqc2} converges uniformly and absolutely on the closed unit disk.

(ii) $\Rightarrow$ (iii). Obvious.

(iii)  $\Rightarrow$ (i).

(IP (i).) Let $x,y \in D$, and assume there exists $T \in \bh$ such that
$Tx = yT$. Let $\vare > 0$, and choose a free polynomial $p$ such that
$\| p(x) - F(x) \| < \vare$ and $\| p(y) - F(y) \| < \vare$.
Then
\begin{eqnarray*}
\| T F(x) - F(y) T\|  &\=& 
\| T F(x) - T p(x) + p(y)T - F(y) T\| \\
&\leq& 2 \| T \| \vare .
\end{eqnarray*}
As $\vare $ is arbitrary, we conclude that $TF(x) = F(y)T$.

(IP (ii).) 
Suppose $(x_n)$ is a bounded sequence in $D$, and assume it is infinite.
(The argument for finite sequences is similar).
Let $z$ be the diagonal $d$-tuple with entries $x_1, x_2 , \dots$,
and let $s : \h \to \h^\i$ be such that
$ y = s^{-1} z s$ is in $D$.
For each fixed $n$, choose a sequence $p_k$ of free polynomials that approximate
$F$ on $\{ y, x_n \}$.
Then
\beq
F \Big( s^{-1}
\begin{bmatrix} x_1 & 0& \cdots \\
0 & x_2 & \cdots\\
\vdots & \vdots & \ddots \end{bmatrix}
s \Big)
&\=&
\lim_{k \to \i} p_k
\Big( s^{-1}
\begin{bmatrix} x_1 & 0& \cdots \\
0 & x_2 & \cdots\\
\vdots & \vdots & \ddots \end{bmatrix}
s \Big) \\
&\=&
s^{-1}
\lim_{k \to \i} 
\begin{bmatrix}p_k  ( x_1) & 0& \cdots \\
0 & p_k ( x_2) & \cdots\\
\vdots & \vdots & \ddots \end{bmatrix} 
s .
\eeq
The $n^{\rm th}$ diagonal entry of the right hand side is $F(x_n)$; so we conclude as $n$ is arbitrary that
\[
F \Big( s^{-1}
\begin{bmatrix} x_1 & 0& \cdots \\
0 & x_2 & \cdots\\
\vdots & \vdots & \ddots \end{bmatrix}
s \Big)
\=
s^{-1}
\begin{bmatrix} F  ( x_1) & 0& \cdots \\
0 & F( x_2) & \cdots\\
\vdots & \vdots & \ddots \end{bmatrix} 
s .
\]

(SSOC). Suppose $x_n $ in $D$ converges to $x$ in $D$ in the SOT. As before, by taking direct sums,
we can approximate $F$ by free polynomials uniformly on countable bounded subsets of $D$.
So for any vector $v$, and any $\vare >0$, we choose a free polynomial $p$ so that
$\|[ F(x_n) - p(x_n) ] v \| < \vare/3$, and choose $N$ so that for $n \geq N$, we have
$\| [p(x) - p(x_n) ]  v \| < \vare/3$. Then 
$\| [F(x) - F(x_n) ]  v \| < \vare$ for all $n \geq N$.
\ep

\vs
 In particular, we
get the following consequence, which says that bounded IP functions leave closed algebras invariant.

In \cite[Thm. 7.7]{amif} it is shown
that  for general nc-functions $f$, it need not be true that $f(x)$ is in the algebra generated by $x$.

\begin{cor}
Assume that $D$ is balanced and closed with respect to countable direct sums, and that $F: D \to \bh$ is SSOC and IP.
Then, for each $x \in D$, the operator $F(x)$ is in the closed unital algebra generated by $x^1, \dots, x^d$.
\end{cor}

\section{Free IP functions}

Recall the definition of the sets $\gdel$ in \eqref{eqa1}; the topology they generate is called the free topology on
$\md$.
\bd
\label{defa33}
A free holomorphic function on a free open set $D \subseteq \md$ is an nc-function that,
in the free topology, is locally bounded.
\ed
 Free holomorphic functions are a class of nc-functions studied by the authors
in \cite{amfree, amfreePick}. In particular, it was shown that there was a representation theorem 
for nc-functions that are bounded by $1$ on $\gdel$.
\bt \label{thmq1}
\cite[Thm. 8.1]{amfree}
Let $\delta$ be an $I$-by-$J$ matrix of free polynomials, and let $f$ be an nc-function on $\gdel$ that is bounded by $1$.
There exists an auxiliary Hilbert space $\L$ and an isometry
\[
  \begin{bmatrix}\alpha&B\\C&D\end{bmatrix} \ : \C \oplus \L^{I} \to \C \oplus \L^{J}
\]
so that for $x \in \gdel \cap B(\K)^d$, 
\be
\label{eqa44}
f(x) \= \alpha I_\K  +  (I_\K \otimes B) (\d(x) \ot I_\L)  \big[ I_\K \otimes I_{\L^J} - (I_\K \otimes D)  (\d(x) \ot I_\L) \big]^{-1}
(I_\K \otimes C).
\ee
\et
Obviously one can define a function  on $\gds$ using the right-hand side of \eqref{eqa44},
replacing the finite dimensional space $\K$ by the infinite dimensional space $\h$.
The following theorem gives sufficient conditions on a function to arise this way.

\bt
\label{thma1}
Let $\delta$ be an $I$-by-$J$ matrix of free polynomials, and assume that 
$\gds$ is connected and contains $0$. Let $F : \gds \to \bho$ be sequentially strong operator continuous. Then the following are equivalent:

(i) The function $F$ is intertwining preserving.

(ii) For each $t > 1$, the function $F$ is uniformly approximable by free polynomials on
$G_{t \d}^\sharp$.

(iii) There exists $\alpha \in \D$ such that if $\Phi = \phi_\alpha \circ F$, then $\Phi^\flat$ is a free holomorphic
function on $\gdel$ that is bounded by $1$ in norm, and that maps $0$ to $0$.

(iv) There exists an auxiliary Hilbert space $\L$ and an isometry
\[
  \begin{bmatrix}\alpha&B\\C&D\end{bmatrix} \ : \C \oplus \L^{I} \to \C \oplus \L^{J}
\]
so that for $x \in \gds$, 
\be
\label{eqa4}
F(x) \= \alpha I_\h  +  (I_\h \otimes B) (\d(x) \ot I_\L)  \big[ I_\h \otimes I_{\L^J} - (I_\h \otimes D)  (\d(x) \ot I_\L) \big]^{-1}
(I_\h \otimes C).
\ee
\et

{\sc Proof of Theorem~\eqref{thma1}:}
(i) $\Rightarrow$ (iii). This follows from Proposition~\eqref{propb1}.

(iii) $\Rightarrow$ (iv).
By Theorem~\eqref{thmq1}, we get such a representation for all $x \in \gdel$.
The series on the right-hand side of \eqref{eqa4} that one gets by expanding the Neumann series of
\be
\label{eql1}
\big[ I_\h \otimes I_{\L^J} - (I_\h \otimes D)  (\d(x) \ot I_\L) \big]^{-1}
\ee
converges absolutely on 
$\gds$; let us denote this limit by $H(x)$.
By Theorem~\eqref{thmc1}, since $H$ is a limit of free polynomials, it is IP and SSOC.
Moreover, as $ \begin{bmatrix}\alpha&B\\C&D\end{bmatrix}$ is an isometry, we get by direct calculation
that
\beq
\lefteqn{ I_\h - H^*(x) H(x) \=} \\
&&
(I_\h \otimes C^*)
 \big[ I_\h \otimes I_{\L^J} -  (\d(x)^* \ot I_\L) (I_\h \otimes D^*) \big]^{-1}
\\
&&
\big[ I_\h \otimes I_{\L^J} -   \d(x)^* \d(x) \ot I_\L \big]\\
&&
 \big[ I_\h \otimes I_{\L^J} - (I_\h \otimes D)  (\d(x) \ot I_\L) \big]^{-1}
(I_\h \otimes C) \\
&& \quad \geq 0.
\eeq
Indeed, to see the last equality without being deluged by tensors, let us write
\[
H(x) \= \dot \alpha  + \dot B \dot \d [ I - \dot D \dot \d]^{-1} \dot C,
\]
where the dots denote appropriate tensors in \eqref{eqa4}.
Then
\beq
I - H^* H &\=& I - \dot \alpha^* \dot \alpha - \dot \alpha^*  \dot B \dot \d [ I - \dot D \dot \d]^{-1} \dot C 
 -   \dot C^*   [ I -  \dot \d^* \dot D^* ]^{-1} \dot \d^* \dot B^* \dot \alpha \\
&& - \dot C^*   [ I -  \dot \d^* \dot D^* ]^{-1}\dot \d^* \dot B^* \dot B \dot \d [ I - \dot D \dot \d]^{-1} \dot C\\
&=&  \dot C^* \dot C + \dot C^* \dot D \dot \d [ I - \dot D \dot \d]^{-1} \dot C + 
\dot C^*  [ I -  \dot \d^* \dot D^* ]^{-1} \dot \d^*  \dot D^* \dot C \\
 &&  - \dot C^*   [ I -   \dot \d^* \dot D^*]^{-1} \dot \d^* [ I -  \dot D^* \dot D ] \dot \d [ I - \dot D \dot \d]^{-1} \dot C \\
 &=&
 \dot C^*   [ I -   \dot \d^* \dot D^*]^{-1}  \left\{  [ I -   \dot \d^* \dot D^*]  [ I - \dot D \dot \d]  \right.  \\
 &&\quad \left.
 +  [ I -   \dot \d^* \dot D^*] \dot D \dot \d + \dot \d^* \dot D^*  [ I - \dot D \dot \d] \right. \\
 &&\quad \left. 
  - \dot \d^*  [ I -  \dot D^* \dot D ] \dot \d \right\}
  [ I - \dot D \dot \d]^{-1} \dot C \\
  &=&
  \dot C^*   [ I -   \dot \d^* \dot D^*]^{-1} [  I -   \dot \d^* \dot \d]  [ I - \dot D \dot \d]^{-1} \dot C .
\eeq

Therefore $\| H(x) \| \leq 1 $ for all $x \in \gds$.
 Let $\Delta(x) = H(x) - F(x)$.
Then $\Delta$ is a bounded IP SSOC  Fr\'echet holomorphic function on  $\gds$ that vanishes on
$\gds \cap \Md = \gdel.$ There is a balanced neighborhood $U$ of $0$ in $\gds$ .
By Theorem~\eqref{thmc1}, $\Delta$ has a power series expansion
$\Delta(x) = \sum P_k(x)$, and each $P_k$ vanishes on $U \cap \md$.
This means each $P_k$ vanishes on a
neighborhood of zero in every $\mnd$, and hence must be zero.
Therefore $\Delta $ is identically zero on $U$. 
By analytic continuation, $\Delta$ is identically zero on all of $\gds$,
and therefore \eqref{eqa4} holds.


(iv) $\Rightarrow$ (ii). This follows because the Neumann series obtained by expanding
\eqref{eql1}
has the $k^{\rm th}$ term bounded by $\| \delta(x) \|^k$.
Therefore it converges uniformly and absolutely on $G_{t \d}^\sharp$ for every $t > 1$.

(ii) $\Rightarrow$ (i). Repeat the argument of (iii) $\Rightarrow$ (i)
of Theorem~\eqref{thmc1}. 
\ep


In the notation of the theorem, let $F^\flat = \phi_{-\alpha} \circ \Phi^\flat$. Then the proof of (iii) $\Rightarrow$ (iv) shows that $F$ and $F^\flat$ determine each other uniquely. So we get:
\bt
\label{thme1}
Let $\delta$ be an $I$-by-$J$ matrix of free polynomials, and assume that 
$\gds$ is connected and contains $0$. Then every bounded free holomorphic function on
$\gdel$ has a unique extension to an IP SSOC function on $\gds$.
\et

\bibliography{../../references}

\end{document}